\documentclass{amsart}
\newcommand{\RN}[1]{%
  \textup{\uppercase\expandafter{\romannumeral#1}}%
}
\usepackage{amsmath}
\usepackage{amssymb}
\usepackage{bm}
\usepackage{enumerate}
\usepackage[hidelinks]{hyperref}

\usepackage{mathtools}
\usepackage[capitalise]{cleveref}
\numberwithin{equation}{section}

\makeatletter
\def\section{\@startsection{section}{1}%
  \z@{.5\linespacing\@plus.7\linespacing}{-.5em}%
  {\normalfont\bfseries}}
\makeatother

\crefname{equation}{}{}
\newtheorem{theorem}{Theorem}[section]

\theoremstyle{definition}

\theoremstyle{remark}
\newtheorem{remark}[theorem]{Remark}

\usepackage{verbatim} 

\usepackage{tikz}
\usepackage{tikz-cd}

\usepackage{enumitem}
\usepackage{epigraph}

\begin{document}
\title{The Wheel Conditions and K-theoretic Hall Algebras}
\author{Yu Zhao}
\address{Beijing Institute of Technology, Beijing, China}
\email{zy199402@live.com} 

\begin{abstract}
In this note, we give a geometric realization of the wheel conditions initiated by Feigin-Odesskii \cite{odesskii1997elliptic,feigin2001quantized} through the $\mathrm{K}$-theoretic Hall algebra on $\mathbb{A}^2$.
\end{abstract}

\maketitle
\section{Motivation}
Schiffmann--Vasserot \cite{schiffmann2013} and Feigin--Tsymbaliuk \cite{feigin2011} constructed the Fock space representation of the quantum toroidal algebra $U_{q,t}(\ddot{\mathfrak{gl}}_{1})$ on the $\mathrm{K}$-theory of Hilbert schemes of points on $\mathbb{A}^{2}$, generalizing the action of Nakajima and Grojnowski on cohomology groups. In \cite{schiffmann2013}, Schiffmann--Vasserot introduced the $\mathrm{K}$-theoretic Hall algebra, an associative algebra on the Grothendieck group of the moduli stack of zero-dimensional coherent sheaves on $\mathbb{A}^{2}$. They proved that the $\mathrm{K}$-theoretic Hall algebra is isomorphic to $U_{q,t}(\ddot{\mathfrak{gl}}_{1})^+$, which acts on the $\mathrm{K}$-theory of Hilbert schemes of points on $\mathbb{A}^{2}$ through Hecke correspondences. In \cite{feigin2011}, Feigin--Tsymbaliuk considered a homomorphism from the positive part $U_{q,t}(\ddot{\mathfrak{gl}}_{1})^+$ to a shuffle algebra $\mathrm{SH}$ with a set of conditions called ``wheel conditions'', and constructed the representation via the shuffle presentation. The definition of the $\mathrm{K}$-theoretic Hall algebra was generalized by Minets \cite{minets2018cohomological} for Higgs sheaves over a curve, and by the author \cite{zhao2019k} and Kapranov--Vasserot \cite{kapranov2019cohomological} for any projective surface $S$. Porta and Sala \cite{porta2019categorification} categorified these Hall algebras.

The purpose of this note is to compare the two presentations of the quantum toroidal algebra $U_{q,t}(\ddot{\mathfrak{gl}}_{1})$ and in particular, to give a geometric interpretation of the wheel conditions in the shuffle algebra presentation of $U_{q,t}(\ddot{\mathfrak{gl}}_{1})$, which can be naturally generalized to the $\mathrm{K}$-theoretic Hall algebra of any quiver $Q$. 

\begin{remark}
  After the first version of the paper was posted on arXiv, the approach of this paper has been extensively generalized by Negu\c{t} \cite{negut2023,NSO2025} to give a presentation of the positive part of the quantum toroidal algebra associated to any quiver $Q$ through the shuffle algebra with wheel conditions. Recently, it was also studied by Gubarevich \cite{Gubarevich} for the wheel conditions of the Hall induction for cotangent representations.
\end{remark}

\section{Main Results}
Let us first recall the construction of Schiffmann--Vasserot's $\mathrm{K}$-theoretic Hall algebra when $S=\mathbb{A}^{2}$. For any $n\in \mathbb{Z}_{>0}$, let $\mathrm{Coh}_{n}$ be the moduli stack of length $n$ coherent sheaves on $\mathbb{A}^{2}$, which is represented by the quotient stack $[\mathrm{Comm}_{n}/\mathrm{GL}_{n}]$ where 
 $$\mathrm{Comm}_{n}=\{(X,Y)\in \mathrm{End}(\mathbb{C}^n)\oplus \mathrm{End}(\mathbb{C}^n)|[X,Y]=0\},$$
is the variety of two commuting matrices, and the action of $\mathrm{GL}_{n}$ on $\mathrm{Comm}_{n}$ is given by conjugation. The scaling action 
\begin{equation*}
  \mathbb{G}_m\times \mathbb{G}_m\times \mathrm{Comm}_n \to \mathrm{Comm}_n, \quad (t_1,t_2,(X,Y))\mapsto (t_1 X, t_2 Y)
\end{equation*}
induces an action of $\mathbb{G}_m\times \mathbb{G}_m$ on $\mathrm{Coh}_n$. For any two non-negative integers $n,m$, let $\mathrm{Corr}_{n,m}$ be the moduli stack of short exact sequences
$$0\to \mathcal{E}_{n}\to \mathcal{E}_{n+m}\to \mathcal{E}_{m}\to 0$$
where $\mathcal{E}_{n}\in \mathrm{Coh}_{n}$, $\mathcal{E}_{m}\in \mathrm{Coh}_{m}$ and
$\mathcal{E}_{n+m}\in \mathrm{Coh}_{n+m}$. Let 
\begin{equation*}
  p:\mathrm{Corr}_{n,m}\to \mathrm{Coh}_{n}\times \mathrm{Coh}_{m}, \quad q:\mathrm{Corr}_{n,m}\to \mathrm{Coh}_{n+m}
\end{equation*}
be the morphisms defined by forgetting and remembering the middle term $\mathcal{E}_{n+m}$, respectively. Schiffmann--Vasserot constructed a virtual pullback morphism
$$p^{!}:\mathrm{K}^{\mathbb{G}_m\times \mathbb{G}_m}(\mathrm{Coh}_{n}\times \mathrm{Coh}_{m})\to
\mathrm{K}^{\mathbb{G}_m\times \mathbb{G}_m}(\mathrm{Corr}_{n,m})$$ 
such that 
$$q_{*}p^{!}:\mathrm{K}^{\mathbb{G}_m\times \mathbb{G}_m}(\mathrm{Coh}_{n}\times \mathrm{Coh}_{m})\to \mathrm{K}^{\mathbb{G}_m\times \mathbb{G}_m}(\mathrm{Coh}_{n+m})$$
induces an associative algebra structure on $\mathrm{KHA}:=\bigoplus_{n=1}^{\infty}\mathrm{K}^{\mathbb{G}_m\times \mathbb{G}_m}(\mathrm{Coh}_{n})$. 

The shuffle algebra $\mathrm{SH}$ is a graded algebra on 
\begin{equation*}
  \mathrm{SH}:=\bigoplus_{n=0}^{\infty}\mathbb{Z}[q_{1}^{\pm 1},q_{2}^{\pm 1}][z_{1}^{\pm 1},\ldots, z_{n}^{\pm 1}]^{\sigma_{n}},
\end{equation*}
where $\sigma_n$ is the symmetric group on $n$ letters. The multiplication is given by the shuffle product, which we will omit but refer to in \cite{feigin2011,schiffmann2013,negut2014} for details. The homomorphism from $\mathrm{KHA}$ to $\mathrm{SH}$ is constructed through the localization theorem in equivariant K-theory \cite{thomason1992formule}. Let $T_{n}\subset \mathrm{GL}(n)$ be the maximal torus, which consists of diagonal matrices. Then we have 
\begin{equation*}
  \mathrm{K}^{\mathbb{G}_m\times \mathbb{G}_m}(\mathrm{Coh}_{n})=(\mathrm{K}^{\mathbb{G}_m\times \mathbb{G}_m\times T_{n}}(\mathrm{Comm}_{n}))^{\sigma_{n}},
\end{equation*}
as $\sigma_{n}$ is the Weyl group of $\mathrm{GL}_{n}$. The fixed locus of $\mathbb{G}_m\times \mathbb{G}_m\times T_n$ on $\mathrm{Comm}_{n}$ is the zero matrix $(0,0)$. Hence the morphism 
$$r^{*}_{0}\circ \psi_{*}:\mathrm{K}^{T_{n}\times \mathbb{G}_{m}\times
  \mathbb{G}_{m} }(\mathrm{Comm}_{n})\to \mathrm{K}^{T_{n}\times \mathbb{G}_{m}\times \mathbb{G}_{m}}(\cdot)=\mathbb{Z}[q_{1}^{\pm 1},q_{2}^{\pm
  1}][z_{1}^{\pm 1},\ldots, z_{n}^{\pm 1}]$$
  is an isomorphism after the localization, where $r_0$ is the inclusion of the point $(0,0)$ in $\mathrm{Comm}_n$ and $\psi$ is the inclusion morphism from $\mathrm{Comm}_n$ to $\mathrm{End}(\mathbb{C}^n)\oplus \mathrm{End}(\mathbb{C}^n)$. By considering the direct sum over all $n$, we have an algebra morphism 
  \begin{equation*}
    \mathrm{KHA}\to \mathrm{SH}:=\bigoplus_{n=0}^{\infty}\mathbb{Z}[q_{1}^{\pm 1},q_{2}^{\pm 1}][z_{1}^{\pm 1},\ldots, z_{n}^{\pm 1}]^{\sigma_{n}}.
  \end{equation*}
  This is an algebra homomorphism. The main result of this paper is the following theorem, which gives a set of wheel conditions for the image of $r^{*}_{0}\circ p_{*}$:
  
\begin{theorem}
    \label{wheel:comm}
   For any three distinct numbers
    $i,j,k\in \{1,\ldots,n\}$, the image of
    $$r^{*}_{0}\circ p_{*}:\mathrm{K}^{T_{n}\times \mathbb{G}_{m}\times
  \mathbb{G}_{m} }(\mathrm{Comm}_{n})\to \mathrm{K}^{T_{n}\times \mathbb{G}_{m}\times \mathbb{G}_{m}}(\cdot)=\mathbb{Z}[q_{1}^{\pm 1},q_{2}^{\pm
  1}][z_{1}^{\pm 1},\ldots, z_{n}^{\pm 1}].$$
is contained in the ideal $(z_{k}-q_{1}z_{j},z_{j}-q_{2}z_{i})$. For the direct sum, the image of the algebra morphism from $\mathrm{KHA}$ to $\mathrm{SH}$ is contained in the intersection of all such ideals for any distinct $i,j,k$.
\end{theorem}

\begin{proof}
  Let $E_{ij}$ be the matrix whose $(i,j)$-th entry is $1$ and
  all other entries are $0$. For three distinct numbers $i,j,k\in
  \{1,\ldots,n\}$, we define
  $$V_{ijk}^{1}=\{(c_{1}E_{ij},0)|c_{1}\in \mathbb{C}\},V_{ijk}^{2}=\{(0,c_{2}E_{jk})|c_{2}\in \mathbb{C}\},V_{ijk}=\{(c_{1}E_{ij},c_{2}E_{jk})|c_{1},c_{2}\in \mathbb{C}\}.$$
 Let
  $$i_{V}:V_{ijk}\to gl_{n}\oplus gl_{n}$$ be the
  inclusion morphism and
  $v_{0}:\mathit{Spec}(\mathbb{C})\to V_{ijk}$
  be the inclusion of the point $(0,0)$. Then $v_{0}^{*}\circ
  i_{V}^{*}=r_{0}^{*}$. Using the following Cartesian diagram:
  \begin{equation}
    \label{diagram:wheel}
    \begin{tikzcd}
      V_{ijk}^{1}\cup V_{ijk}^{2} \ar{r}{p'} \ar{d} & V_{ijk} \ar{d}{i_{V}} \\
      \mathrm{Comm}_{n} \ar{r}{p} & gl_{n}\oplus gl_{n},
    \end{tikzcd}
  \end{equation}
we consider the refined Gysin morphism 
\begin{align*}
  i_V^!:\mathrm{K}^{\mathbb{G}_m\times \mathbb{G}_m\times T_n}(\mathrm{Comm}_n)\to \mathrm{K}^{\mathbb{G}_m\times \mathbb{G}_m\times T_n}(V_{ijk}^1\cup V_{ijk}^2) \\
  [\mathcal{F}]\to \sum_{i=0}^{\infty}(-1)^i[Tor^i_{\mathcal{O}_{gl_n\oplus gl_n}}(\mathcal{O}_{V_{ijk}},\mathcal{F})].
\end{align*}
By the definition of the refined Gysin morphism and the Cartesian diagram \cref{diagram:wheel}, we have $ r_{0}^{*}\circ p_{*}= v_{0}^{*}\circ i_{V}^{*}\circ p_{*}=v_{0}^{*}\circ p_{*}' \circ i_{V}^{!}$. Thus the image of $r_{0}^{*}\circ p_{*}$ is contained in the image of
  $v_{0}^{*}\circ p_{*}'$. The image of $p_{*}'$ is the ideal
  generated by two classes $[\mathcal{O}_{V_{ijk}^{1}}]$ and
  $[\mathcal{O}_{V_{ijk}^{2}}]$, while
  $$v_{0}^{*}([\mathcal{O}_{V_{ijk}^{1}}])=1-q_{2}^{-1}\frac{z_{k}}{z_{j}}, \ v_{0}^{*}([\mathcal{O}_{V_{ijk}^{2}}])=1-q_{1}^{-1}\frac{z_{j}}{z_{i}}.$$
  Hence the image of $v_{0}^{*}\circ p_{*}'$ is contained in the ideal  $(z_{k}-q_{2}z_{j},z_{j}-q_{1}z_{i})$.
    
\end{proof}

\section{Acknowledgements} The author would like to thank Andrei Negut for the enjoyable time spent talking with him and for all the knowledge about shuffle algebras learned from him. The author would also like to thank Svetlana Makarova for much help in learning algebraic $\mathrm{K}$-theory and Francesco Sala for many useful suggestions.

\end{document}